\documentclass[12pt, reqno]{amsart}

\usepackage[T2A]{fontenc}
\usepackage[cp1251]{inputenc}
\usepackage[english]{babel}
\usepackage{amsfonts,amssymb, amsthm} 

\theoremstyle{plain}

\newtheorem{condition}{Condition}
\newtheorem{theorem}{Theorem}
\newtheorem{proposition}{Proposition}
\newtheorem{remark}{Remark}
\newtheorem{col}{Colorary}

\newcommand\const{\operatorname{const}}
\newcommand\Imm{\operatorname{Im}}
\newcommand\Ree{\operatorname{Re}}

\newcommand{\Id}{{I}}
\newcommand{\Emb}{{\mathcal J}}
\newcommand{\Hil}{{\mathfrak H}}
\newcommand{\R}{{\mathbb R}}
\newcommand{\Complex}{{\mathbb C}}
\newcommand{\Dom}{{\mathcal D}}
\newcommand{\BOp}{{\mathcal{L} }}
\newcommand{\T}{{\mathcal{T} }}

\author{Nikita\,V.~Artamonov}
\email{nikita.artamonov@gmail.com}

\title[On exponential decay rate of semigroup]{On exponential decay rate of semigroup associated with 
second order linear differential equation in Hilbert space with strong damping operator}


\thanks{This paper is supported by the Russian Foundation of Basic Research (project No
11-01-00790)}

\begin{document}

\begin{abstract}
We obtain  estimate of the exponential decay rate of
semigroup associated with second order linear differential
equation $u''+Du'+Au=0$ in Hilbert space. We assume that $A$
is a selfadjoint positive definite operator,  $D$ is an accretive
sectorial operator and $\Ree D\geq\delta A$, $\delta>0$.
We obtain a location of the spectrum of a pencil associated with
linear differential equation.

\end{abstract}

\maketitle


\textbf{1.}  In Hilbert space $H, (\cdot,\cdot)$ we consider a second order linear differential equation
 \begin{equation}\label{SeconfOrderDiffEqn}
  u''(t)+Du'(t)+Au(t)=0,
 \end{equation}
here $u(t)$ is a vector-value function on semi-axis $\R_+=[0,+\infty)$. Many evolution equations arising in mechanics
 can be reduced to the equation \eqref{SeconfOrderDiffEqn} in an appropriate
space (see, for example,  \cite{CramerLatushkin03, Huang97, JacobTrunk09}). In this case $A$ represents 
potential energy and $D$ represents dissipation ($D$ is a damping operator). We will assume that
\begin{condition}\label{ConditionA}
$A$ is a selfadjoit positive definite operator with dense domain $\Dom(A)$. Let
 \[
     a_0=\inf_{x\in\Dom(A),\|x\|=1}(A x,x)=
     \inf_{x\in\Dom(A^{1/2}),\|x\|=1}(A^{1/2} x,A^{1/2}x)>0.
 \]
\end{condition}
By $H_s$ we denote a collection of Hilbert spaces generated by $A^{1/2}$, i.e. for $s\geq0$ the space $H_s$ is the
domain $\Dom(A^{s/2})$ endowed with the norm $\|x\|_s=\|A^{s/2}x\|$, for $s<0$ the space $H_s$ is the completion of $H$
with respect to the norm $\|\cdot\|_s$. By definition $H_0=H,H_1=\Dom(A^{1/2}),H_2=\Dom(A)$ and 
$H_s\hookrightarrow H_r$ for $s>r$. Since $|(x,y)|=|(A^{-s/2}x,A^{s/2}y)|\leq\|x\|_{-s}\cdot\|y\|_s$  for $s>0$ and for all 
$x\in H, y\in H_s$, then the sesquilinear  form $(x,y)$ can be extended by continuity to a sesquilinear form 
$(x,y)_{-s,s}$ on $H_{-s}\times H_s$. Therefore we can regard the space $H_{-s}$ as a dual of $H_s$ ($H_{-s}=H_s^*$), 
the duality is determined by  the sesquilinear form $(x,y)_{-s,s}$ (duality with respect to the pivot space $H$). By $\BOp(X,Y)$
we will denote a space of bounded operators acting from a space $X$ into a space $Y$. Operator $A$ can be regarded
as a bounded operator acting in the collection of Hilbert spaces: $A\in\BOp(H_s,H_{s-2})$ $\forall s\in\R$.

Following \cite{JacobTrunk09} we will assume that
\begin{condition}\label{ConditionD}
$D\in\BOp(H_1,H_{-1})$ is an accretive sectorial operator, i.e. $\Ree(Dx,x)_{-1,1}\geq0$ and 
$|\Imm(Dx,x)_{-1,1}|\leq\nu \Ree(Dx,x)_{-1,1}$ for all $x\in H_1$  and some $\nu>0$.
\end{condition}
Denote (infimum with respect to $x\in H_1,x\ne0$)
 \[
    \alpha=\inf\frac{\Ree(Dx,x)_{-1,1}}{\|x\|^2_{-1}},\;
    \beta=\inf\frac{\Ree(Dx,x)_{-1,1}}{\|x\|^2},\;
    \delta=\inf\frac{\Ree(Dx,x)_{-1,1}}{\|x\|^2_{1}}.
 \]
Inequality $\|x\|^2_1\geq a_0\|x\|^2\geq a_0^2\|x\|^2_{-1}$ ($\forall x\in H_1$)  implies, that $$\alpha\geq a_0\beta\geq a_0^2\delta\geq0.$$
By $\|D\|=\sup_{x\in H_1,\|x\|_1=1}\|Dx\|_{-1}$ we denote a norm of the operator  $D$. Note, that the operator $D\in\BOp(H_1,H_{-1})$ 
is accretive (sectorial) in the sense of the condition \ref{ConditionD} iff the operator $A^{-1/2}DA^{-1/2}\in\BOp(H)$ 
is accretive (sectorial). 

With the linear differential equation \eqref{SeconfOrderDiffEqn} we associate a quadratic operator
pencil \cite{HrynivShkalikov04,JacobTrunk09}
 \[
    L(\lambda)=\lambda^2 \Emb+\lambda D+A,
 \]
here $\Emb:H_1\hookrightarrow H_{-1}$ is an embedding operator, $\lambda\in\Complex$ is a spectral parameter.
We regard the pencil as an operator-function $L(\lambda)\in\BOp(H_1,H_{-1})$. As usual  one can define a resolvent set
  \[
     \rho(L)=\{\lambda\in\Complex\;:\;\exists L^{-1}(\lambda)\in\BOp(H_{-1},H_1)\}
  \]
and a spectrum $\sigma(L)=\Complex\backslash\rho(L)$ of the pencil $L(\lambda)$.

The second order differential equation \eqref{SeconfOrderDiffEqn} can be linearized as a first order
differential equation \cite{HrynivShkalikov04,Huang97,JacobTrunk09}
\begin{equation}\label{FirstOrderSystem}
 w'(t)=\T w(t), \quad w(t)=\begin{pmatrix} u' & u \end{pmatrix}^\top
\end{equation} 
in "energy" space $\Hil=H\times H_1$ with the matrix operator 
 \[
    \T=\begin{pmatrix}
     -D & -A \\ \Id & 0
    \end{pmatrix},\;
    \Dom(\T)=\left\{\begin{pmatrix} w_1 \\ w_2 \end{pmatrix}\in H_1\times H_1\,:\,
    D w_1+A w_2\in H
    \right\}
 \]
Since $\Ree(\T w,w)_{\Hil}=-\Ree(Dw_1,w_1)_{-1,1}\leq0$ and $0\in\rho(\T)$, then ($-\T$) is a maximal
accretive operator and, therefore, $\T$ is a generator of $C_0$-semigroup  $\exp(t\T)$ of contractions
\cite{EngelNagel2000}. In \cite{HrynivShkalikov04} it is shown that if the conditions \ref{ConditionA} and
\ref{ConditionD} hold and $\beta>0$, then the operator $\T$ is a generator of exponentially decaying semigroup. 
In \cite{HrynivShkalikov99}, in particular, it is  proved that, under the conditions \ref{ConditionA},
\ref{ConditionD} and $\delta>0$, the operator $\T$ is a generator of analytic semigroup.

In papers \cite{BatkaiEngel04, Huang97} for the case $\beta>0$ was obtained results on the exponential decay rate of the
semigroup  $\exp(t\T)$ and on location of the spectrum of the pencil $L(\lambda)$. In paper \cite{JacobTrunk09} was obtained
results on analyticity of the semigroup $\exp(t\T)$ and on the location of the spectrum of the pencil $L(\lambda)$.
In present paper using another technique for the cases $\delta>0$ we obtain an estimate for the exponential decay rate of
the semigroup generated by $\T$.

\textbf{2.} In the space $\Hil$ with respect to the given inner product the operator $(-\T)$ is neither uniformly accretive nor sectorial.
For $\theta\geq0$  introduce a collection of sesquilinear forms
 \begin{multline*}
 [w,v]_\theta=(w_1,v_1)+\theta(w_1,v_1)_{-1}+(w_2,v_2)_1+\theta (w_2,v_2)+\theta(Dw_2,Dv_2)_{-1}+\\
 \theta(Dw_2,v_1)_{-1}+\theta(w_1,Dv_2)_{-1},\quad w=\begin{pmatrix} w_1 \\ w_2 \end{pmatrix},\,
  v=\begin{pmatrix} v_1 \\ v_2 \end{pmatrix}\in\Hil,
 \end{multline*}
here $(\cdot,\cdot)_{s}=(A^{s/2}\cdot,A^{s/2}\cdot)$ is an inner product in the space $H_s$. 
Since
 \[
   |w|^2_\theta=[w,w]_\theta=\|w_1\|^2+\|w_2\|^2_1+\theta\|w_2\|^2+\|w_1+Dw_2\|^2_{-1},
 \]
then $[\cdot,\cdot]_\theta$ is an inner product in $\Hil$ topologically equivalent to
the given one. Obviously $[\cdot,\cdot]_0=(\cdot,\cdot)_{\Hil}$.
\begin{proposition}\label{QuaraticFormEstimate}
Let the conditions \ref{ConditionA} and \ref{ConditionD} are satisfied and $\delta>0$. 
Then for arbitrary $\theta>0$ and $0\leq b\leq\sqrt{\theta}$ for all 
$w=(w_1,w_2)^\top\in\Dom(\T)$ the following inequalities
 \[
    \Ree[\T w,w]_\theta\leq -\omega_\theta|w|^2_\theta,\quad 
    \bigl|\Imm[\T w,w]_\theta\bigr|\leq M_{\theta,b}\bigl|\Ree[\T w,w]_\theta\bigr|+b|w|^2_\theta,
 \] 
hold, where
 \begin{gather*}
   \frac{1}{\omega_\theta} =\frac{1}{\beta}+\frac{\|D\|^2}{2\delta}+\frac12\left(\frac{\theta}{\alpha}+\frac{1}{\theta\delta}\right)+
   \frac12\sqrt{\left(\frac{\theta}{\alpha}+\frac{1}{\theta\delta}+\frac{\|D\|^2}{\delta}\right)^2-\frac{4}{\alpha\delta}}>0, \\
   M_{\theta,b} = \nu+\frac{2}{\delta(b+\sqrt{b^2+4\theta})}+\frac{\sqrt{\theta}-b}{\beta}.
 \end{gather*}
\end{proposition}
The resolvent set of the operator $\T$ is non-empty, therefore 
\begin{col}
Under the conditions of the proposition \ref{QuaraticFormEstimate}
 \[
    \sigma(\T)\subset\{\lambda\in\Complex\;:\; \Ree\lambda\leq -\omega_\theta, 
     |\Imm\lambda|\leq M_{\theta,b}|\Ree\lambda|+b\}.
 \]
\end{col}
Putting $b=0$ we obtain
\begin{col}
Under the conditions of the proposition \ref{QuaraticFormEstimate}
for all $\theta>0$
 \[
    \sigma(\T)\subset\left\{\lambda\in\Complex\;:\;\Ree\lambda\leq -\omega_\theta, 
    |\Imm\lambda|\leq  \left(\nu+\frac{1}{\delta\sqrt{\theta}}+\frac{\sqrt{\theta}}{\beta}\right)
    |\Ree\lambda|\right\}
 \]
\end{col}
It's easy to prove \cite{JacobTrunk09}, that $\sigma(L)=\sigma(\T)$ .
\begin{theorem}\label{MainResult}
Let the conditions \ref{ConditionA} and \ref{ConditionD} are satisfied and $\delta>0$.
Then the operator $\T$ is a generator of the (analytic) semigroup $\exp(t\T)$ in the space 
$\Hil$ with exponential decay rate 
\[
    \omega =\left(\frac{1}{\beta}+\frac{2}{\sqrt{\alpha\delta}}+\frac{\|D\|^2}{2\delta}+
   \sqrt{\frac{4\|D\|^2}{\delta\sqrt{\alpha\delta}}+\frac{\|D\|^4}{\delta^4}}\right)^{-1}>0,
 \]
i.e. for all $t\geq0$ the inequality $\|\exp(t\T)\|_\Hil\leq \const\cdot\exp(-\omega t)$ holds. 
For all $b\geq0$
 \[
     \sigma(L)=\sigma(\T)\subset\{\lambda\in\mathbb{C}\;|\; \Ree\lambda\leq -\omega,
    |\Imm\lambda|\leq M_b|\Ree\lambda|+b\}
 \]
where  $M_b=\min_{\theta\geq b^2}M_{\theta,b}$. 
\end{theorem}
\begin{col}
Under the conditions of the theorem \ref{MainResult} for all $(u_1,u_0)^\top\in\Dom(\T)$ there exists a
unique solution $u(t)$ of the Cauchy problem for the differential equation
\eqref{SeconfOrderDiffEqn} with initial conditions $u(0)=u_0$, $u'(0)=u_1$
and 
 \[
    \|u(t)\|^2_1+\|u'(t)\|^2\leq\const\cdot\exp(-2\omega t)(\|u_0\|^2_1+\|u_1\|^2).
 \]
\end{col}
Putting $b=0$ we have
\begin{col}
Under the conditions of the theorem \ref{MainResult}
 \[
     \sigma(L)=\sigma(\T)\subset\{\lambda\in\mathbb{C}\;|\; \Ree\lambda\leq -\omega,
    |\Imm\lambda|\leq \left(\nu+\frac{2}{\sqrt{\delta\beta}}\right)|\Ree\lambda|\}.
 \]
\end{col}
\begin{remark}
In \cite{JacobTrunk09} under the condition of the theorem \ref{MainResult} was obtained the following location
of the pencil's spectrum
 \begin{align*}
     \sigma(L)=\sigma(\T)&\subset\{\lambda\in\Complex\;|\; \Ree\lambda\leq0,\;
    |\Imm\lambda|\leq \nu|\Ree\lambda|+\delta^{-1}\}\\
     \sigma(L)=\sigma(\T)&\subset\left\{\lambda\in\Complex\;|\;\delta\leq\frac{|\Ree\lambda|}{a_0^{-2}+|\lambda|^{-2}}\right\}
 \end{align*}
\end{remark}

\textbf{3.} With the operator $\T$ in the space $\Hil$ endowed with the inner
product $[\cdot,\cdot]_\theta$ we can associate a linearization of the pencil
$L(\lambda)$ under the form $\mathbf{L}(\lambda)=\lambda \mathbf{Q}-\mathbf{T}$,
where
 \[ 
   \mathbf{Q}=\begin{pmatrix} \Id+\theta A^{-1} & \theta A^{-1}D \\ \theta D^* A^{-1} & A+\theta\Id+\theta D^* A^{-1}D \end{pmatrix}
   \quad
   \mathbf{T}=\begin{pmatrix} -D & -A-\theta\Id \\ A+\theta\Id &   -\theta D^* \end{pmatrix}.
 \]
The linearization  $\mathbf{L}(\lambda)$ can be regarded as an operator-function
$\mathbf{L}(\lambda)\in\BOp(H_1\times H_1,H_{-1}\times H_{-1})$.
 
\textbf{Acknowledgement} The author thanks prof. A.A.~Shkalikov and prof. C.~Trunk 
for fruitful discussions.

\end{document}